\documentclass{article}
\usepackage{amssymb}

\usepackage{amsfonts,wasysym}
\usepackage{graphicx}

\newtheorem{theorem}{Theorem}[section]

\newtheorem{definition}[theorem]{Definition}

\newtheorem{example}[theorem]{Example}

\newtheorem{remark}[theorem]{Remark}

\newtheorem{conjecture}[theorem]{Conjecture}

\newtheorem{problem}[theorem]{Problem}

\begin{document}
\vspace*{1.5pc}

\large
{\bf\Large \begin{center}
Braids and branched coverings of dimension three
\end{center}}
\medskip
\centerline{\Large 
J. Scott Carter}

\medskip
\centerline{\Large 
Department of Mathematics, University of South Alabama}

\medskip
\centerline{and} 
\medskip

\centerline{\Large 
Seiichi Kamada}

\medskip
\centerline{\Large 
Department of Mathematics, Hiroshima University}

\bigskip
\medskip

%%%%%%%%%%%%%%%%%%%%%
\section{Introduction}
%%%%%%%%%%%%%%%%%%%%%

This is on a part of our work in progress, which was 
introduced at the conference ``Intelligence of Low-dimensional Topology'' held in RIMS in May, 2012. 
The purpose of our research is to understand branched coverings and $m$-dimensional braids which are generalizations of  classical braids. Here we discuss chart descriptions of branched coverings and braids in dimension $m=2$ first, and then those for which $m=3$. 

We work in the PL category (\cite{Hud1969, RS1972}).  
Let $S^m$ denote the $m$-sphere, and let $M^m$ denote a closed oriented $m$-manifold.  

\section{Preliminaries}

We start by giving some definitions and theorems on branched coverings.  

\begin{definition}{\rm 
A PL map $f: M^m \to S^m$ is a {\it branched covering} ({\it map}) if 
there exists an $(m-2)$-subcomplex $L$ of $S^m$ such that the restriction 
$\underline f : M^m \setminus f^{-1}(L) \to S^m \setminus L$ is a 
covering map.  
}\end{definition}  

We denote the covering degree by $d$.  We call $f$ a {\it $d$-fold} branched covering.  

We assume that $L$ is minimum, i.e., 
$\forall y \in L$, 
$\#(f^{-1}(y)) < d$.   Then we call $L$ the {\it branch set} of $f$.   

\begin{definition}{\rm 
A $d$-fold branched covering  $f$ is {\it simple} if $\forall y \in L$, $\#(f^{-1}(y)) = d-1$.  
}\end{definition}

\begin{remark}{\rm 
(1) A branched covering is defined in general as follows (cf. \cite{BE1979, BE1984}): 
A PL map between manifolds is called {\it proper} if the inverse image of the boundary is the boundary.  
A proper PL map between manifolds $f : M^m \to N^m$ is called a branched covering if it is finite-to-one and open.  

(2) A branched covering $f: M \to N$ is {\it primitive} if $f_\ast : \pi_1(M) \to \pi_1(N)$ is surjective.  
It is often assumed that a branched covering is primitive. 
}\end{remark}

Note  that $M^m$ is closed, oriented and \underline{connected} in what follows in this section.  

\begin{theorem}[J.W. Alexander \cite{Al1920}]
For any closed oriented and connected $m$-manifold $M^m$, there exists a simple branched covering $f : M^m \to S^m$ for some degree $d$.  
\end{theorem}

\begin{remark}{\rm 
(1) A closed oriented and connected $1$-manifold $M^1$ is homeomorphic to $S^1$.  Thus 
there exists a $1$-fold covering $ f : M^1\to S^1$. 

(2) For any closed oriented and connected $2$-manifold $M^2$, there exists a $2$-fold simple branched covering $f : M^2\to S^2$.  
}\end{remark}

\begin{theorem}[H. M. Hilden \cite{Hi1976}, J. M. Montesinos \cite{Mo1976}] 
For any closed oriented and connected $3$-manifold $M^3$,  there exists a $3$-fold simple branched covering $f : M^3\to S^3$ such that  
 the branch set $L$ is a link (or a knot). 
\end{theorem} 

The following is a conjecture due to Montesinos.  

\begin{conjecture}
For any closed oriented and connected $4$-manifold  $M^4$,  there exists a $4$-fold simple branched covering
$f : M^4\to S^4$ such that $L$ is an embedded surface in $S^4$.    
\end{conjecture} 

Some partial answers to this conjecture are known as follows.

\begin{theorem}[R. Piergallini \cite{Pi1995}] 
For any closed oriented and connected $4$-manifold  $M^4$,  there exists a $4$-fold simple branched covering
$f : M^4\to S^4$ such that $L$ is an immersed surface in $S^4$.     
\end{theorem} 

\begin{theorem}[M. Iori and R. Piergallini \cite{IP2002}] 
For any closed oriented and connected $4$-manifold  $M^4$,  there exists a $5$-fold simple branched covering
$f : M^4\to S^4$ such that $L$ is an embedded surface in $S^4$.  
\end{theorem} 

%%%%%%%%%%%%%%%%%%%%%
\section{Two dimensional case $(m=2)$}
%%%%%%%%%%%%%%%%%%%%%

Let 
$f: M^2 \to S^2$ be a $d$-fold simple branched covering with branch set $L$, and let 
$\underline f : M^2 \setminus f^{-1}(L) \to S^2 \setminus L$ be the associated covering map.  

Take a base point $\ast$ of $S^2 \setminus L$ to consider the fundamental group $\pi_1(S^2 \setminus L, \ast)$. 
The preimage $f^{-1}(\ast)$ of the base point $\ast$ consists of $d$ points of $M^2$.  
Then we have a  {\it monodromy } $ \rho: \pi_1(S^2 \setminus L, \ast) \to S_d$, where the symmetric group $S_d$  on letters $\{1, 2, \dots, d\}$  is identified with the symmetric group on $f^{-1}(\ast)$.   (A monodromy $\rho$ depends on the identification between  $\{1, 2, \dots, d\}$ and $f^{-1}(\ast)$.)   The covering $\underline f $ is determined by the monodromy.  

By the Riemann-Hurwitz formula,  $L$ consists of an even number of points.

\begin{figure}[htb]
\begin{center}
\includegraphics[width=5in]{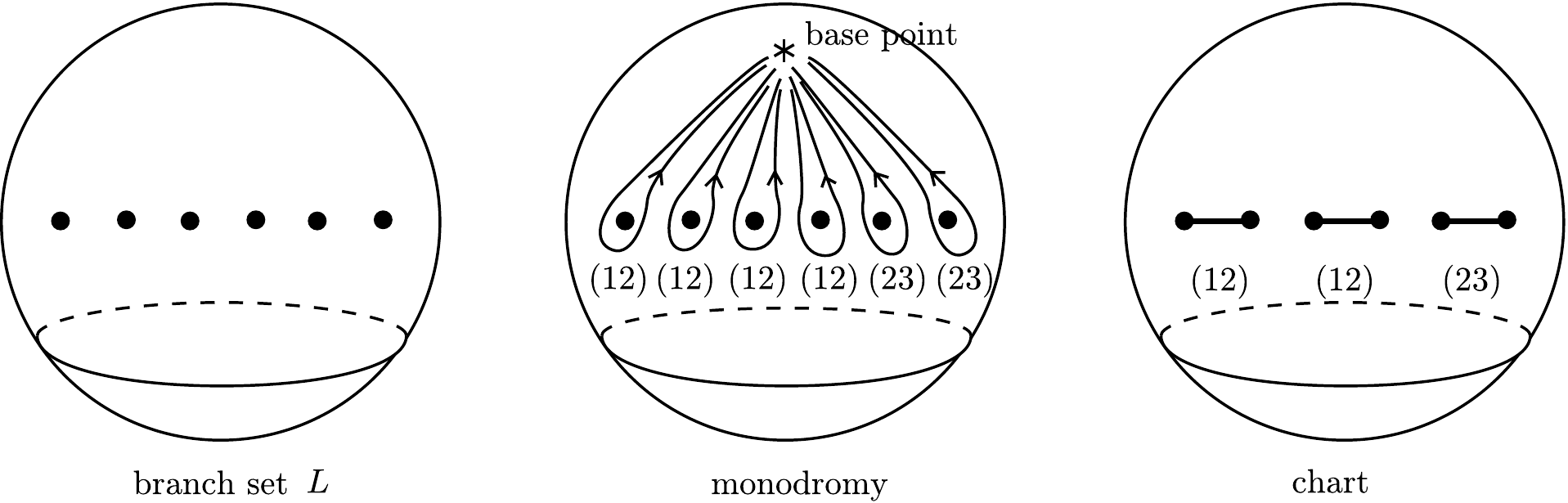}
\end{center}
\caption{A branch set, a monodromy and a chart}
\label{fig:s2_01}
\end{figure}

In Figure~\ref{fig:s2_01}, a branch set, a monodromy, and a chart are depicted.   (A chart description is explained later.)  

\begin{figure}[htb]
\begin{center}
\includegraphics[width=5in]{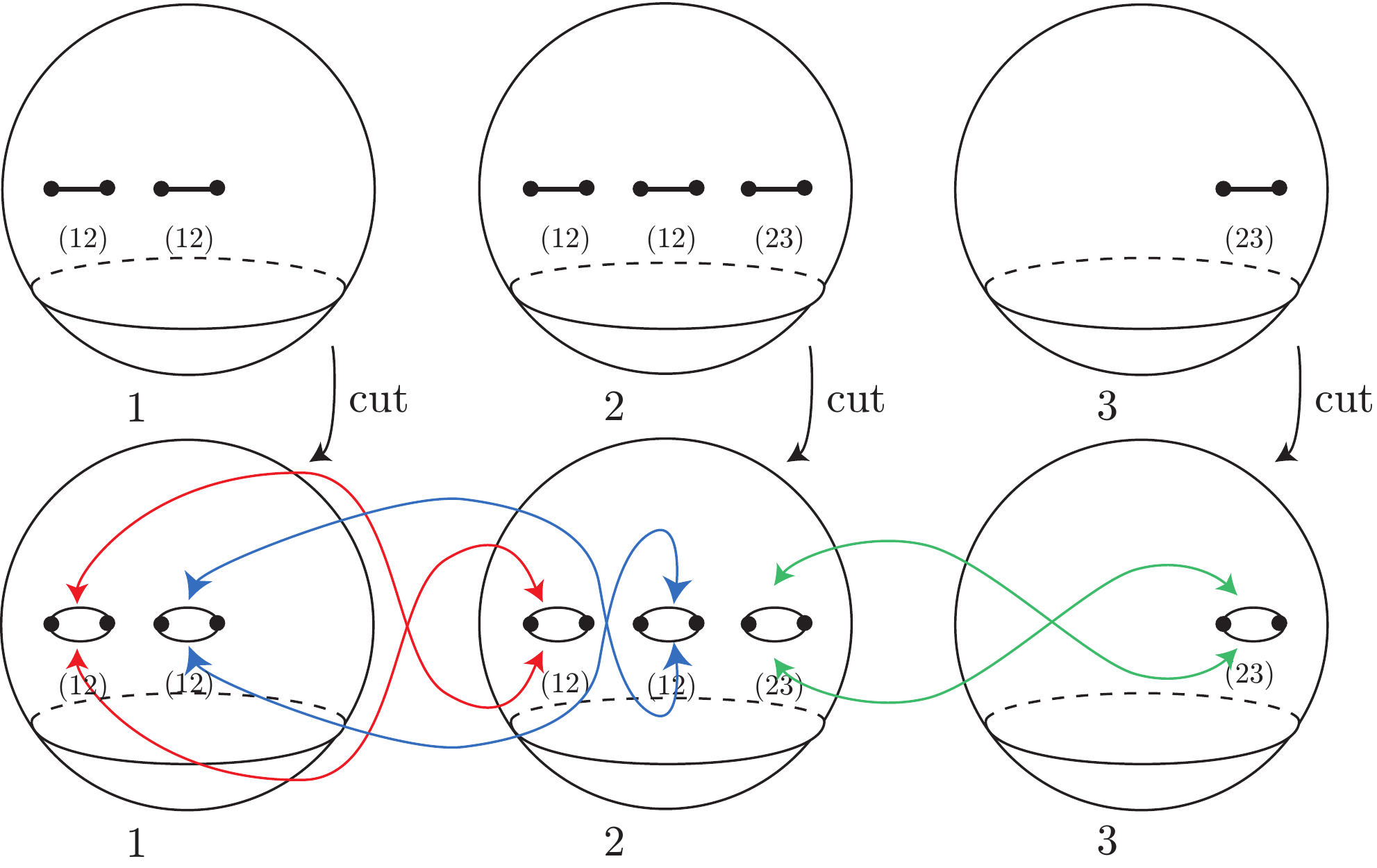}
\end{center}
\caption{How to construct $M^2$}
\label{fig:s2_04}
\end{figure}

When a monodromy is described by a chart, it is easy to 
construct $M^2$. We explain it by using an example.   
Let $\Gamma$ be the chart depicted on the right of Figure~\ref{fig:s2_01}.  
Consider three copies of $S^2$ labeled by $1$, $2$, and $3$, say $S^2_1$, $S^2_2$ and $S^2_3$, respectively.    
On the copy $S^2_1$, draw the edges with label $(12)$ of $\Gamma$, 
on the copy $S^2_2$, draw the edges with label $(12)$ of $\Gamma$ and those with label $(23)$, 
and on the copy $S^2_3$, draw the edges with label $(23)$.  
Cut the three $2$-spheres along these edges, and we obtain three compact surfaces, 
say $M_1$, $M_2$ and $M_3$,  as in the bottom of Figure~\ref{fig:s2_04}.  
The surface $M^2$ is obtained from the union $M_1 \cup M_2 \cup M_3$  by identifying the boundary as follows:  Let $e$ be an edge 
with label $(12)$ on $S^2_1$, and let $e_+$ and $e_-$ be the copies of $e$ in $\partial M_1$.  
Let $e'$ be the corresponding edge on $S^2_2$, and let $e'_+$ and $e'_-$ be the corresponding copies in $\partial M_2$.  Then we identify $e_+$ with $e'_-$, and identify $e_-$ with $e'_+$, respectively.  All boundary edges of $M_1 \cup M_2 \cup M_3$ are identified in this fashion, and we have a closed surface.  This is the desired $M^2$.

The classification of simple branched coverings was studied by 
J. L{\" u}roth \cite{Lur1871}, A. Clebsch \cite{Cl1973},  A. Hurwitz \cite{Hur1891}, and others. 
The classification theorem is stated as follows.  

\begin{theorem}\label{thm:classification1}
Let $f : M^2 \to S^2$ and $f' : {M^2}' \to S^2$ be $d$-fold simple branched coverings with branch sets  $L$ and $L'$, respectively.  We assume that $M^2$ and ${M^2}' $ are connected.  Then 
$ f $ and $ f'$ are equivalent if and only if $\# L = \# L'$.  
\end{theorem}

Hurwitz  \cite{Hur1891} studied branched coverings by using of a system of monodromies of meridian elements of the branch set, called a 
{\it Hurwitz system}, and studied when two systems present the same (up to equivalence) branched coverings.  

A Hurwitz system depends on a system of generating set of $\pi_1(S^2 \setminus L, \ast)$.  
For a generating system depicted in the middle of Figure~\ref{fig:s2_01}, the 
Hurwitz system is  
$$\alpha = ((12), (12), (12), (12), (23), (23)).$$   
Besides a choice of a generating system, a Hurwitz system  depends  on the identification of 
$\{1, 2, \dots, d\}$ and the fiber $f^{-1}(\ast)$.  

Two Hurwitz systems present the same (up to equivalence) braid monodromy if and only if they are 
related by a finite sequence of  {\it Hurwitz moves}  and {\it conjugations}.  
The {\it Hurwitz moves} are 
$$(a_1, \dots, a_k, a_{k+1}, \dots, a_n) \mapsto (a_1, \dots, a_{k+1}, a_{k+1}^{-1} a_k a_{k+1}, \dots, a_n)$$ 
for $k=1, \dots, n-1$ and their inverse moves.   {\it Conjugations} are 
$$(a_1, \dots,  a_n) \mapsto (g^{-1} a_1 g, \dots, g^{-1} a_n g)$$
for $g \in S_d$.   When two Hurwitz systems are related by a finite sequence of  {Hurwitz moves  and conjugations, we say that they are {\it HC-equivalent}.  ($H$ and $C$ stand for Hurwitz and conjugation.) 

Due to Hurwitz  \cite{Hur1891}, the classification theorem is stated as follows.  

\begin{theorem}
Let $f : M^2 \to S^2$ be a $d$-fold simple branched covering.  Assume that $M^2$ is connected.  
Any Hurwitz system of $f$ is HC-equivalent to 
$$((12), \dots, (12), (23), (23), (34), (34), \dots, (d-1,d), (d-1,d)).$$ 
(The number of $(12)$s is a positive even number, and for each $i=2, \ldots, d-1$, a pair of $(i, i+1)$ appears.)
\end{theorem}

In the next section, we will  introduce the notion of a {\it chart}, called a 
{\it permutation chart} or an {\it $S_d$-chart},  that describes a branched covering or its monodromy.  
The chart method helps us to 
construct $M^2$ from a monodromy, and to understand the classification theorem well.  

%%%%%%%%%%%%%%%%%%%%%
\section{Permutation charts or $S_d$-charts $(m=2)$}
%%%%%%%%%%%%%%%%%%%%%

We denote by $\tau_i$ the transposition $(i ~ i+1)$.  The symmetric group  $S_d$ is generated by 
$\tau_1, \dots, \tau_{d-1}$, and has a group presentation 
$$S_d = \left\langle \tau_1, \dots, \tau_{d-1}  \, 
\begin{array}{|ll}
\tau_i \tau_j \tau_i = \tau_j \tau_i \tau_j  \quad & (|i-j| = 1) \\ 
\tau_i \tau_j = \tau_j \tau_i &   (|i-j| > 1) \\ 
\tau_i^2 = e 
\end{array}
\right\rangle. $$

\begin{definition}{\rm 
A {\it permutation chart} of degree $d$ or an 
 {\it $S_d$-chart}  is a labeled graph in $S^2$ such that each edge is labeled in $\{1, \dots, d-1\}$ and each vertex is as in Figure~\ref{fig:schartvert_01}. We call a vertex a {\it black vertex}, a {\it crossing} or a {\it white vertex} if the valency of the vertex is $1$, $4$ or $6$, respectively.   
}\end{definition}

By the correspondence $i  \leftrightarrow  \tau_i = (i ~  i+1) \in S_d$, the labels of a chart are assumed to be 
transpositions in $S_d$ (see Figure~\ref{fig:s2_01}).     Figure~\ref{fig:scharteg_02} is an example of an $S_4$-chart, or a permutation chart of degree~$4$.

\begin{figure}[htb]
\begin{center}
\includegraphics[width=3.5in]{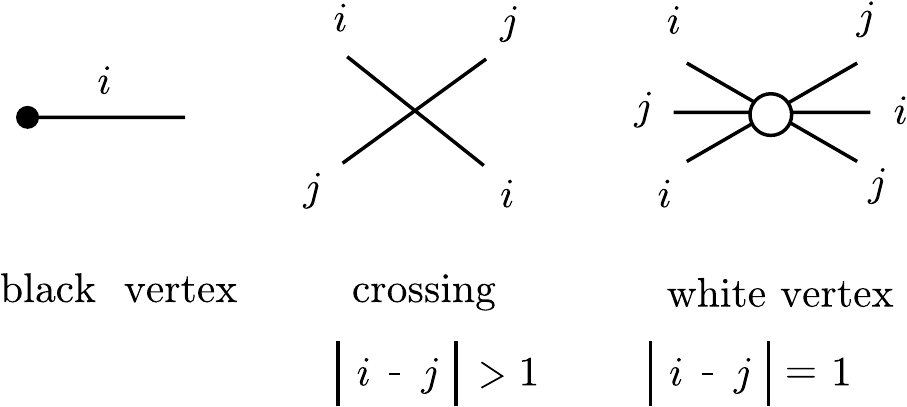}
\end{center}
\caption{Vertices of a $S_d$-chart}
\label{fig:schartvert_01}
\end{figure}

\begin{figure}[htb]
\begin{center}
\includegraphics[width=4in]{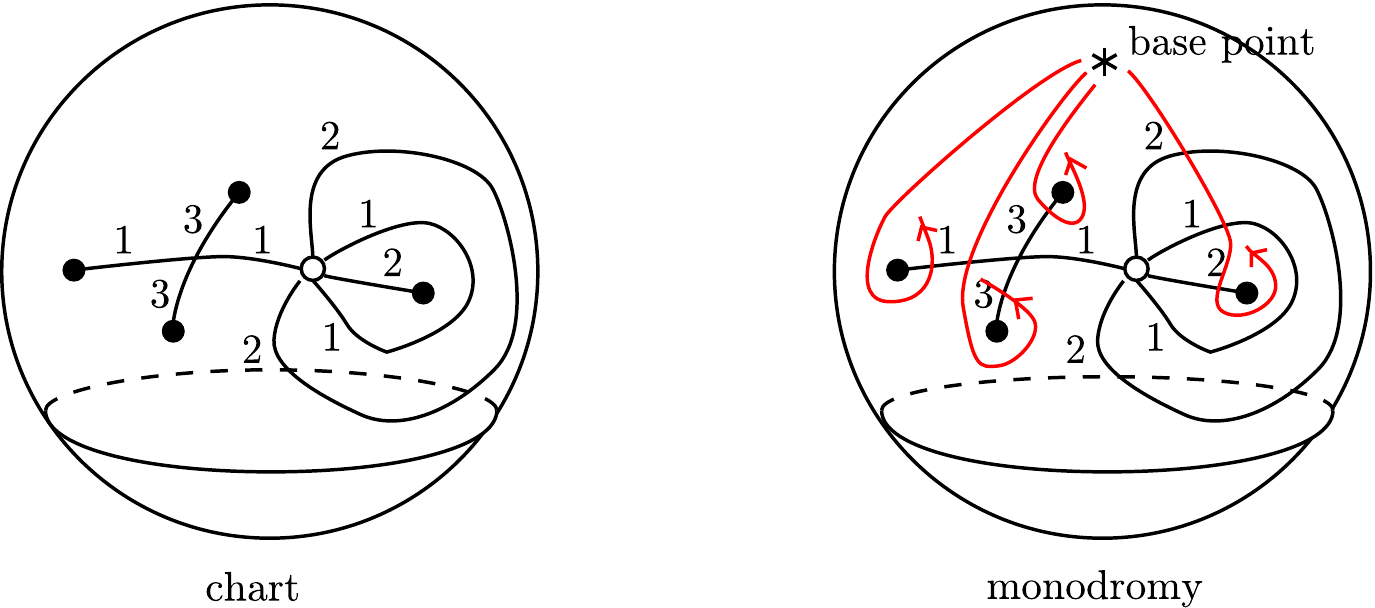}
\end{center}
\caption{A $S_4$-chart $\Gamma$ and the induced monodromy $\rho_\Gamma$}
\label{fig:scharteg_02}
\end{figure}

For a chart $\Gamma$, we consider a monodromy 
$$\rho_\Gamma: \pi_1(S^2 \setminus L) \to S_d,  \quad [\ell] \mapsto [  \textrm{intersection word of $\ell$ w.r.t. $\Gamma$}   ],  $$ 
where $L$ $(= L_\Gamma)$ is the set of black vertices.  An intersection word is a sequence of elements of  $\{1, \dots, d-1\}$, which is regarded  as an element of $S_d$ by the correspondence $i  \leftrightarrow  \tau_i = (i ~  i+1) \in S_d$.

\begin{example}{\rm 
Let $\Gamma$ be an $S_4$-chart depicted in the left of Figure~\ref{fig:scharteg_02}.  
When we take a Hurwitz generating system as in the figure, we have a Hurwitz 
system  $(\tau_1, \, \tau_1 \tau_3 \tau_1, \, \tau_3, \, \tau_2 \tau_1 \tau_2 \tau_1 \tau_2)$.   
It is equal to $(\tau_1, \,  \tau_3, \, \tau_3, \, \tau_1)$.  And it is Hurwitz equivalent to $(\tau_1, \,  \tau_1, \, \tau_3, \, \tau_3)$.  
}\end{example}

\begin{theorem} 
Let  $f : M^2 \to S^2$ be a $d$-fold simple branched covering, and  $\rho_f$ a monodromy of $f$.  
There exists a chart $\Gamma$ such that $\rho_\Gamma = \rho_f$.  
{\rm (We call $\Gamma$ a {\it chart description} of $f$ or $\rho_f$.)}
\end{theorem} 

Local moves on permutation charts illustrated in Figure~\ref{fig:schartmoves} are called {\it chart moves}. 
 (Ignore the orientations on edges.) 
 Two charts are said to be {\it equivalent} or {\it chart move equivalent} if they are related by 
 a finite sequence of chart moves and ambient isotopies of $S^2$.  

\begin{figure}[htb]
\begin{center}
\includegraphics[width=5.5in]{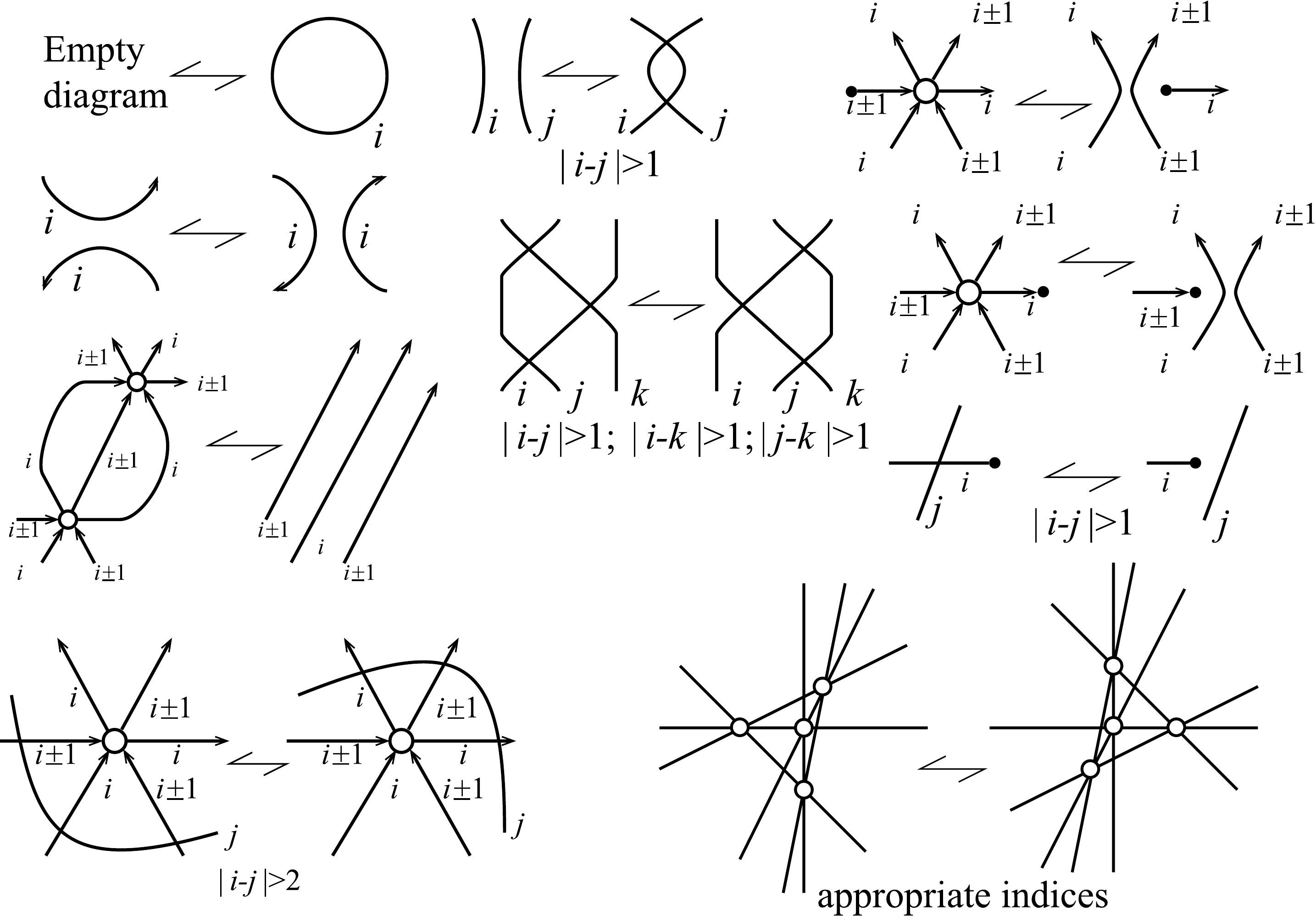}
\end{center}
\caption{Chart moves}
\label{fig:schartmoves}
\end{figure}

\begin{theorem} 
Let $f$ and $f'$ be $d$-fold simple branched covering of $S^2$, and 
let $\Gamma$ and $\Gamma'$ be their chart descriptions.  
$f $ is equivalent to $f'$ if and only if  $\Gamma$ is equivalent to $\Gamma'$.   
\end{theorem} 

Using an example, we explain how to construct $M^2$ from a chart description.  
Let $\Gamma$ be an $S_4$-chart depicted in the top of Figure~\ref{fig:scharteg_03}.   
Consider four copies of $S^2$ labeled by $1$, $2$, $3$ and $4$, say $S^2_1$, $S^2_2$, $S^2_3$  and $S^2_4$, respectively.    
On the copy $S^2_1$, draw the edges with label $1$ of $\Gamma$, 
on the copy $S^2_2$, draw the edges with label $1$ of $\Gamma$ and those with label $2$, 
on the copy $S^2_3$, draw the edges with label $2$ of $\Gamma$ and those with label $3$, 
and on the copy $S^2_4$, draw the edges with label $3$.  
Cut the four $2$-spheres along the edges, and we obtain compact surfaces, 
say $M_1$, $M_2$, $M_3$ and $M_4$,  as in the bottom of Figure~\ref{fig:scharteg_03}.  
The surface $M^2$ is obtained from the union $\cup_{i=1}^{4} M_i$  by identifying the boundary as follows:  Let $e$ be an edge 
with label $1$ on $S^2_1$, and let $e_+$ and $e_-$ be the copies of $e$ in $\partial M_1$.  
Let $e'$ be the corresponding edge on $S^2_2$, and let $e'_+$ and $e'_-$ be the corresponding copies in $\partial M_2$.  Then we identify $e_+$ with $e'_-$, and identify $e_-$ with $e'_+$, respectively.  All boundary edges of $\cup_{i=1}^{4} M_i$ are identified in this fashion, and we have a closed surface.  This is the desired $M^2$.  

At a white vertex, 3 sheets are gathering as in Figure~\ref{fig:strip1}. 

\begin{figure}[htb]
\begin{center}
\includegraphics[width=5in]{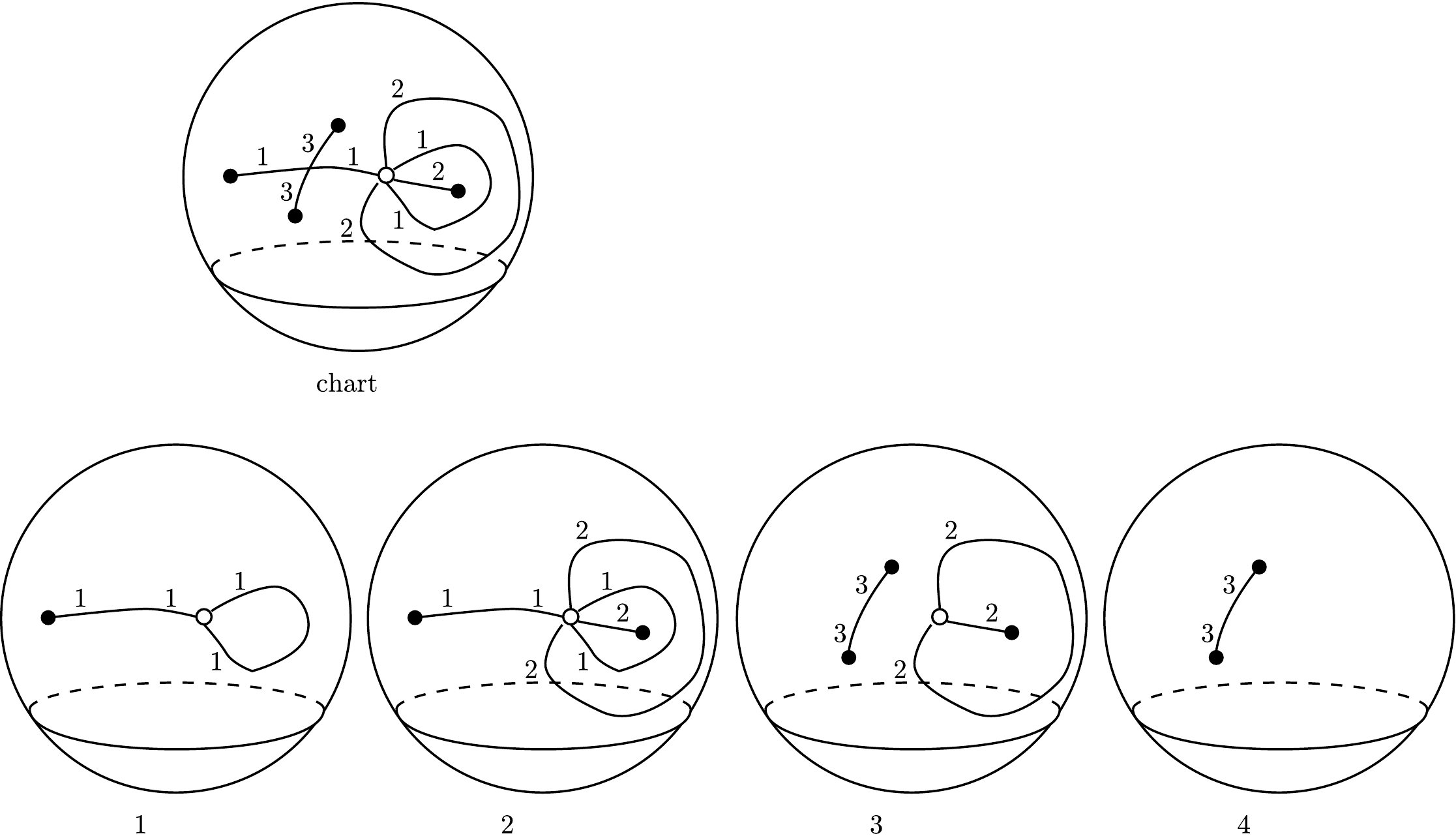}
\end{center} \vspace{-0.1in}
\caption{How to construct $M^2$}
\label{fig:scharteg_03}
\end{figure}

\begin{figure}[htb]
\begin{center}
\includegraphics[width=3.0in]{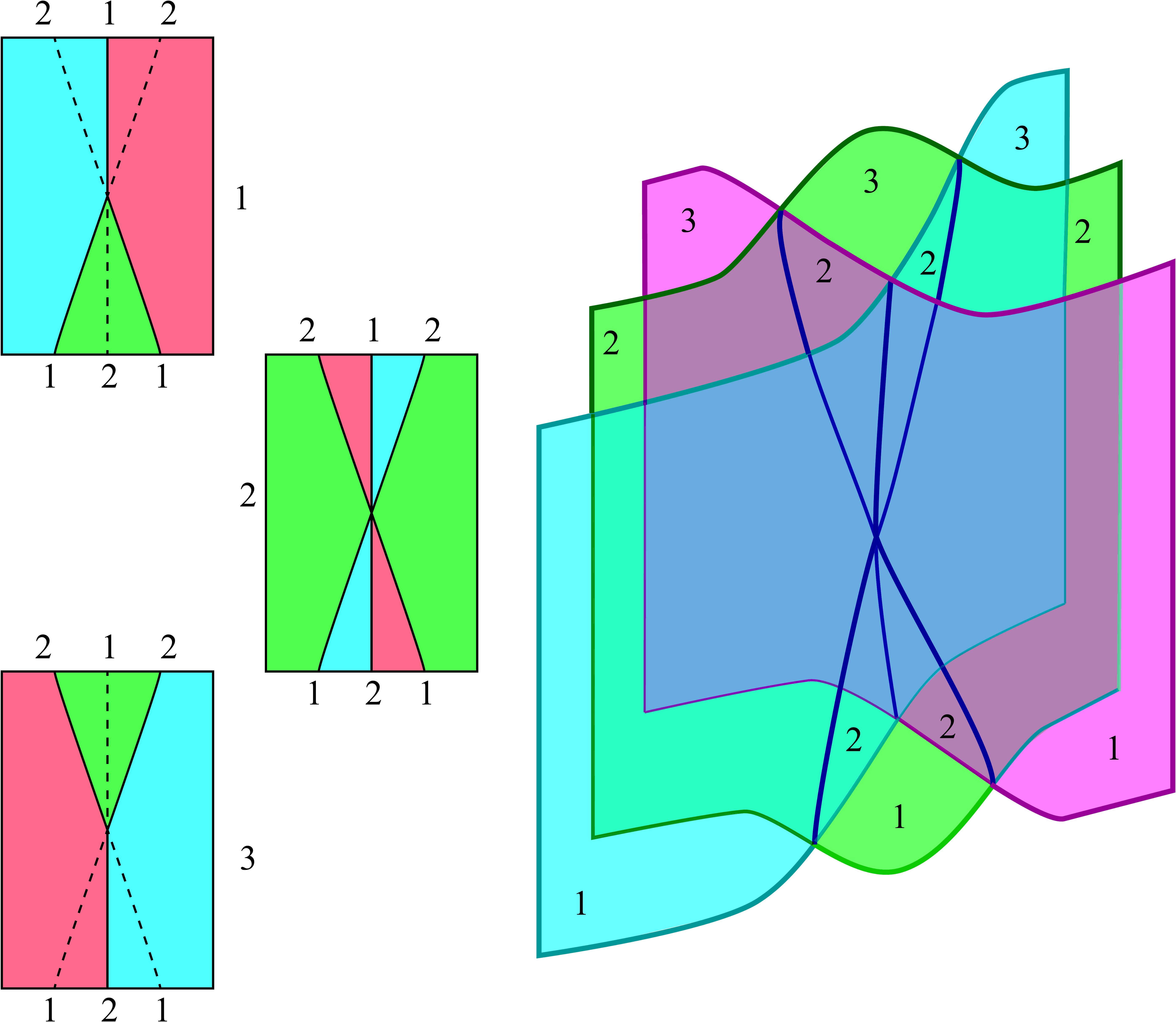}
\end{center}
\caption{Three sheets gather around a white vertex.}
\label{fig:strip1}
\end{figure}

\begin{theorem}
Any chart description of $f : M^2 \to S^2$ with connected $M$ is equivalent to a chart as in Figure~\ref{fig:s2_05}. 
\end{theorem} 

\begin{figure}[htb]
\begin{center}
\includegraphics[width=5in]{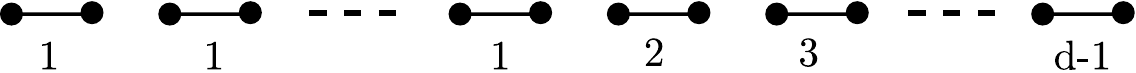}
\end{center}
\caption{A chart in a normal form}
\label{fig:s2_05}
\end{figure}

This theorem is quite easily proved.  
As a corollary of this theorem, we have the classification theorem (Theorem~\ref{thm:classification1}).  

\clearpage 

%%%%%%%%%%%%%%%%%%%%%
\section{Braid charts or $B_d$-charts $(m=2)$}
%%%%%%%%%%%%%%%%%%%%%

Let $\sigma_i$ $(i=1, \dots, d-1)$ be the standard generators of the braid group $B_d$.  Then $B_d$ has a group presentation 
$$B_d = \left\langle \sigma_1, \dots, \sigma_{d-1}  \, 
\begin{array}{|ll}
\sigma_i \sigma_j \sigma_i = \sigma_j \sigma_i \sigma_j  \quad & (|i-j| = 1) \\ 
\sigma_i \sigma_j = \sigma_j \sigma_i &   (|i-j| > 1) 
\end{array}
\right\rangle. $$

\begin{definition}{\rm 
A {\it braid chart} of degree $d$ or a {\it $B_d$-chart}  is a labeled and oriented graph in $S^2$ such that each edge is labeled in $\{1, \dots, d-1\}$ and each vertex is as in Figure~\ref{fig:schartvert_01ori}. We call a vertex a {\it black vertex}, a {\it crossing} or a {\it white vertex} if the valency of the vertex is $1$, $4$ or $6$, respectively.  The arrow at a black vertex in this figure is suppressed since it may either be incoming or outgoing. 
}\end{definition}

\begin{figure}[htb]
\begin{center}
\includegraphics[width=3.5in]{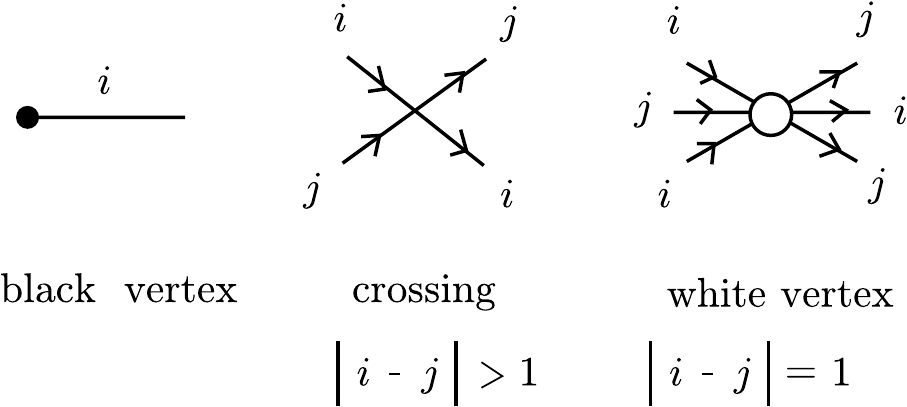}
\end{center}
\caption{Vertices of a $B_d$-chart}
\label{fig:schartvert_01ori}
\end{figure}

By the correspondence $i  \leftrightarrow  \sigma_i = (i ~  i+1) \in B_d$, the labels of a chart are assumed to present 
the standard generators in $B_d$.     Figure~\ref{fig:scharteg_02ori} is an example of a $B_4$-chart, or a braid chart of degree~$4$. 

\begin{figure}[htb]
\begin{center}
\includegraphics[width=4in]{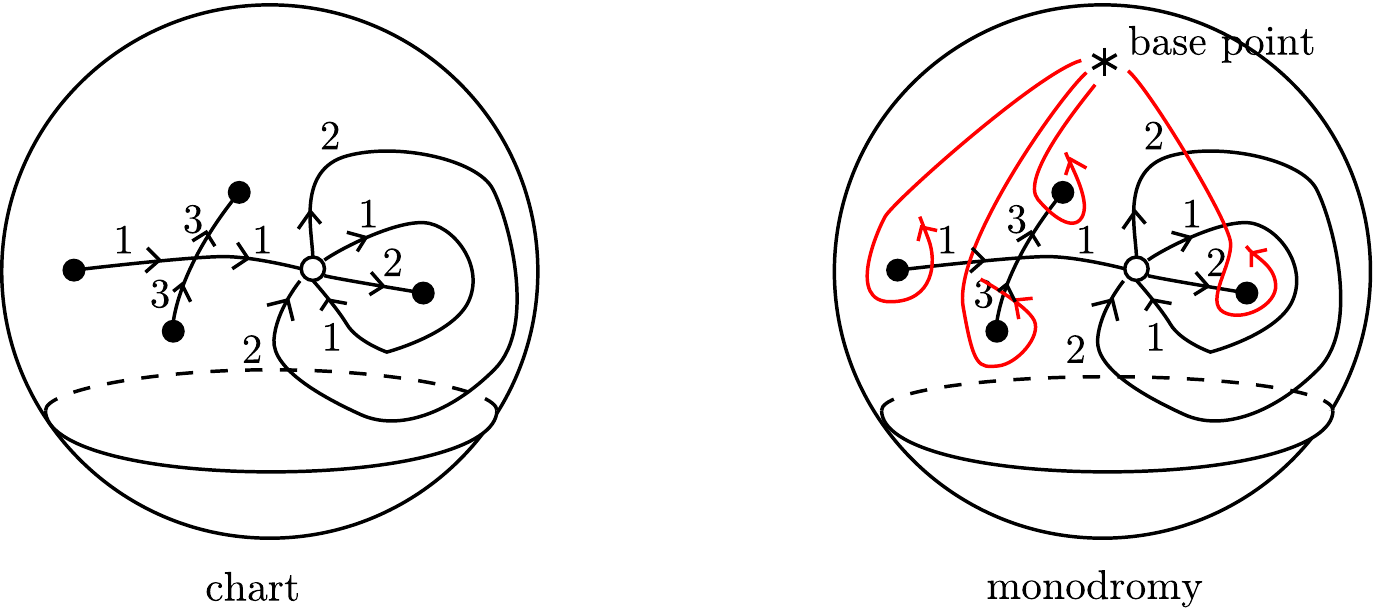}
\end{center}
\caption{A $B_4$-chart $\Gamma$ and the induced monodromy $\rho_\Gamma$}
\label{fig:scharteg_02ori}
\end{figure}

Forgetting orientations of the edges from a braid chart, we obtain a permutation chart. 
Thus we often call a permutation chart an {\it unoriented chart}, and a braid chart an {\it oriented chart}.  

\begin{definition}{\rm 
A permutation chart  is called {\it orientable} if one can give orientations to the edges to make it a braid chart.  
Otherwise it is called {\it nonorientable}. 
}\end{definition}

For a braid chart $\Gamma$ of degree $d$, we consider a monodromy 
$$\rho_\Gamma: \pi_1(S^2 \setminus L) \to B_d,  \quad [\ell] \mapsto [  \textrm{intersection word of $\ell$ w.r.t. $\Gamma$}   ],  $$ 
where $L$ $(= L_\Gamma)$ is the set of black vertices.  An intersection word is a word of  $\{1, \dots, d-1\}$, which is regarded as an element of $B_d$ by the correspondence $i  \leftrightarrow  \sigma_i = (i ~  i+1) \in S_d$.

\begin{example}{\rm 
Let $\Gamma$ be a $B_4$-chart depicted in the left of Figure~\ref{fig:scharteg_02ori}.  
When we take a Hurwitz generating system as in the right of the figure, we have a Hurwitz 
system  
$$(\sigma_1, \, \sigma_1^{-1} \sigma_3 \sigma_1, \, \sigma_3^{-1}, \, \sigma_2^{-1} \sigma_1^{-1} \sigma_2^{-1} \sigma_1 \sigma_2).$$    
It is equal to $(\sigma_1, \,  \sigma_3, \, \sigma_3^{-1}, \, \sigma_1^{-1})$.  And it is Hurwitz equivalent to $(\sigma_1, \,  \sigma_1^{-1}, \, \sigma_3, \, \sigma_3^{-1})$.  
}\end{example}

Let $D^2 \times S^2$ be a tubular neighborhood of a standardly embedded $2$-sphere in $R^4$.  

\begin{definition}{\rm 
A PL embedding $g : M^2 \to D^2 \times S^2 \subset R^4$ is a ({\it simple}) {\it embedded}  {\it  $2$-dimensional braid}, or a {\it surface braid}, of degree $d$  if 
the composition $M^2 \to D^2 \times S^2 \to S^2$ is a $d$-fold (simple) branched covering. 
}\end{definition}

For a (simple or nonsimple) embedded 2-dimensional braid $g : M^2 \to D^2 \times S^2 \subset R^4$ of degree $m$, we can consider a  
 {\it monodromy } $ \rho \, (=\rho_g) : \pi_1(S^2 \setminus L, \ast) \to B_d$, where $L \, (=L_g)$ is the branch set of the branched covering 
$M^2 \to D^2 \times S^2 \to S^2$.  

\begin{theorem}
For any simple embedded $2$-dimensnional braid $g : M^2 \to D^2 \times S^2 \subset R^4$,  there exists a braid chart $\Gamma$ such that $\rho_g = \rho_\Gamma$. 
{\rm ($\Gamma$ is called a {\it chart description} of $g$.)}
\end{theorem} 

Two charts are {\it equivalent} or {\it chart move equivalent} if they are related by a finite sequence of chart moves (Figure~\ref{fig:schartmoves}) and ambient  isotopes of $S^2$.  

\begin{theorem} 
Let $\Gamma$ and $\Gamma'$ be chart descriptions of simple embedded $2$-dimensional braids $g$ and $g'$ of the same degree. 
$g$ and $g'$ are equivalent if and only if  $\Gamma$ is equivalent to $ \Gamma'$.   
\end{theorem} 

Let ${\rm pr}: D^2 \times S^2 \to S^2$ be the projection.  

Let $f: M^2 \to S^2$ be a simple branched covering, and 
$g : M^2 \to D^2 \times S^2$ a simple embedded  2-dimensional braid.  

\begin{definition}{\rm 
If $ {\rm pr} \circ g = f$, 
then we call $g$ an {\it embedded lift}  of $f$, and we say that $f$ is {\it liftable}.  
}\end{definition}

\begin{theorem} 
Any simple branched covering of $S^2$ is liftable.  
\end{theorem} 

\begin{remark}{\rm 
For any simple branched covering, there exists a  chart description that is an orientable permutation chart.   
Not every chart description of a liftable simple branched covering is orientable.  
}\end{remark}

For further topics related to braid charts and $2$-dimensional braids, refer to \cite{CKS2004, CS1998, Kam2002, Kam2012}.  

%%%%%%%%%%%%%%%%%%%%%
\section{Three dimensional case $(m=3)$}
%%%%%%%%%%%%%%%%%%%%%

We recall the theorem due to H. M. Hilden \cite{Hi1976} and J. M. Montesinos \cite{Mo1976} again.  

\begin{theorem}[Hilden  and Montesinos] 
Any closed oriented and connected  $3$-manifold can be represented as a $3$-fold simple branched covering of $S^3$ branched over a link (or a knot).  
\end{theorem} 

Let 
$f: M^3 \to S^3$ be  a $d$-fold simple branched covering of $S^3$ branched along $L$.  
Let $\underline f : M^3 \setminus f^{-1}(L) \to S^3 \setminus L$ be the associated  covering.  
The covering map $\underline f $ is determined by a monodromy 
$ \rho: \pi_1(S^3 \setminus L, \ast) \to S_d$.

\begin{remark}{\rm 
The monodromy $\rho$ sends each meridian to a transposition.  Conversely, 
any homomorphism $ \rho: \pi_1(S^3 \setminus L, \ast) \to S_d$ sending 
each meridian to a transposition is a monodromy of a simple branched covering.  
}\end{remark} 

Figure~\ref{fig:s3_01} is a knot with a monodromy in $S_3$.  In general, 
by 
$(12) \mapsto B={\rm blue}$, 
$(23) \mapsto R={\rm red}$, 
$(13) \mapsto G ={\rm green}$, 
we obtain a link with Fox's 3-coloring that represents a 3-manifold.  See Figure~\ref{fig:scovermove_02}.

\begin{figure}[htb]
\begin{center}
\includegraphics[width=1.5in]{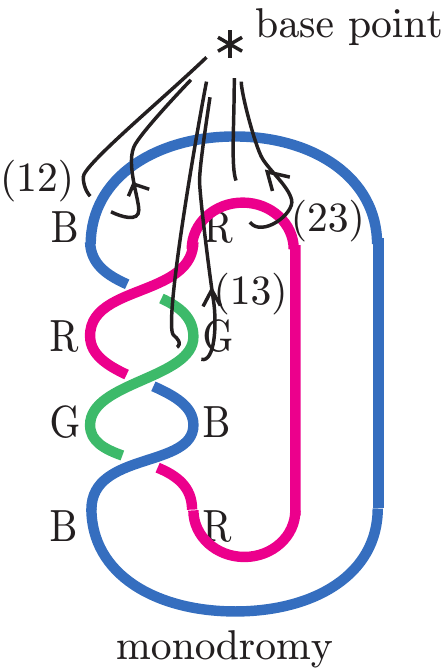}
\end{center}
\caption{A knot with a monodromy in $S_3$}
\label{fig:s3_01}
\end{figure}

\begin{figure}[htb]
\begin{center}
\includegraphics[width=1.2in]{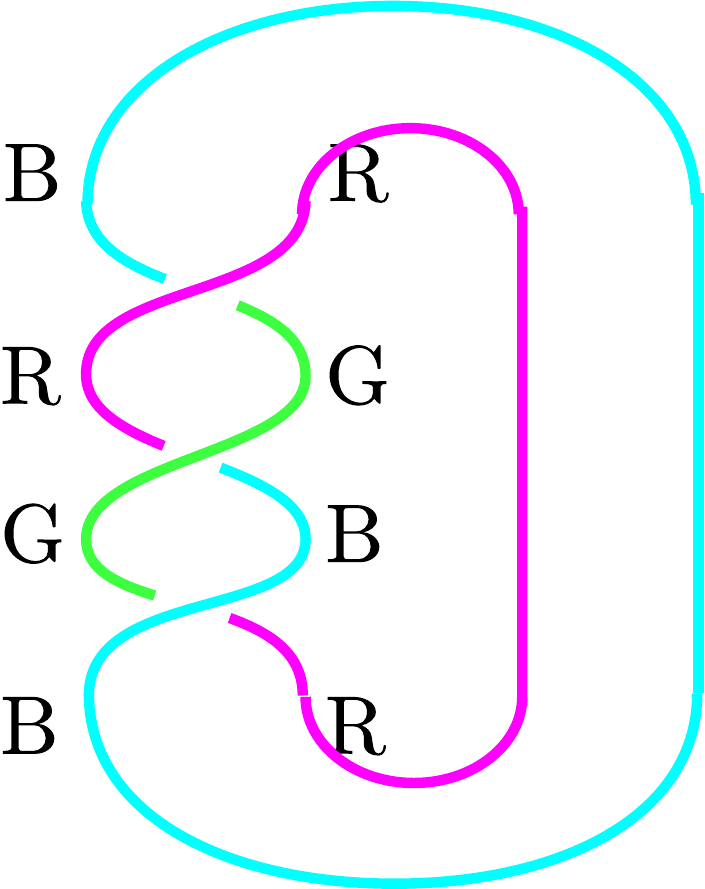}
\end{center}
\caption{A $3$-colored knot}
\label{fig:scovermove_02}
\end{figure}

The  local move depicted in Figure~\ref{fig:scovermove_01} was introduced by Montesinos, that does not change the 3-manifold.  

\begin{figure}[htb]
\begin{center}
\includegraphics[width=2.5in]{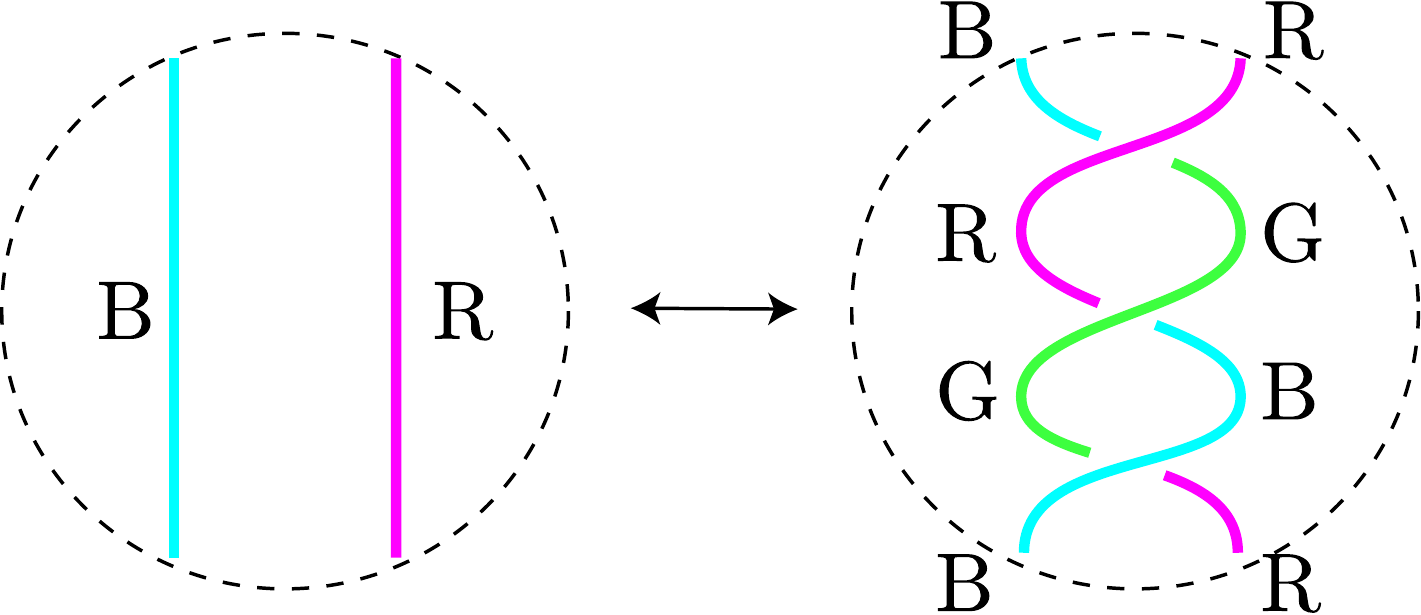}
\end{center}
\caption{A Montesions move}
\label{fig:scovermove_01}
\end{figure}

Applying a Montesions move to the $3$-colored knot in Figure~\ref{fig:scovermove_02}, 
we have a $3$-colored trivial link as in Figure~\ref{fig:scovermove_03}, which represents $S^3$.  
Thus it is a nontrivial representation of $S^3$ as a 3-fold simple branched covering.  

\begin{figure}[htb]
\begin{center}
\includegraphics[width=2.8in]{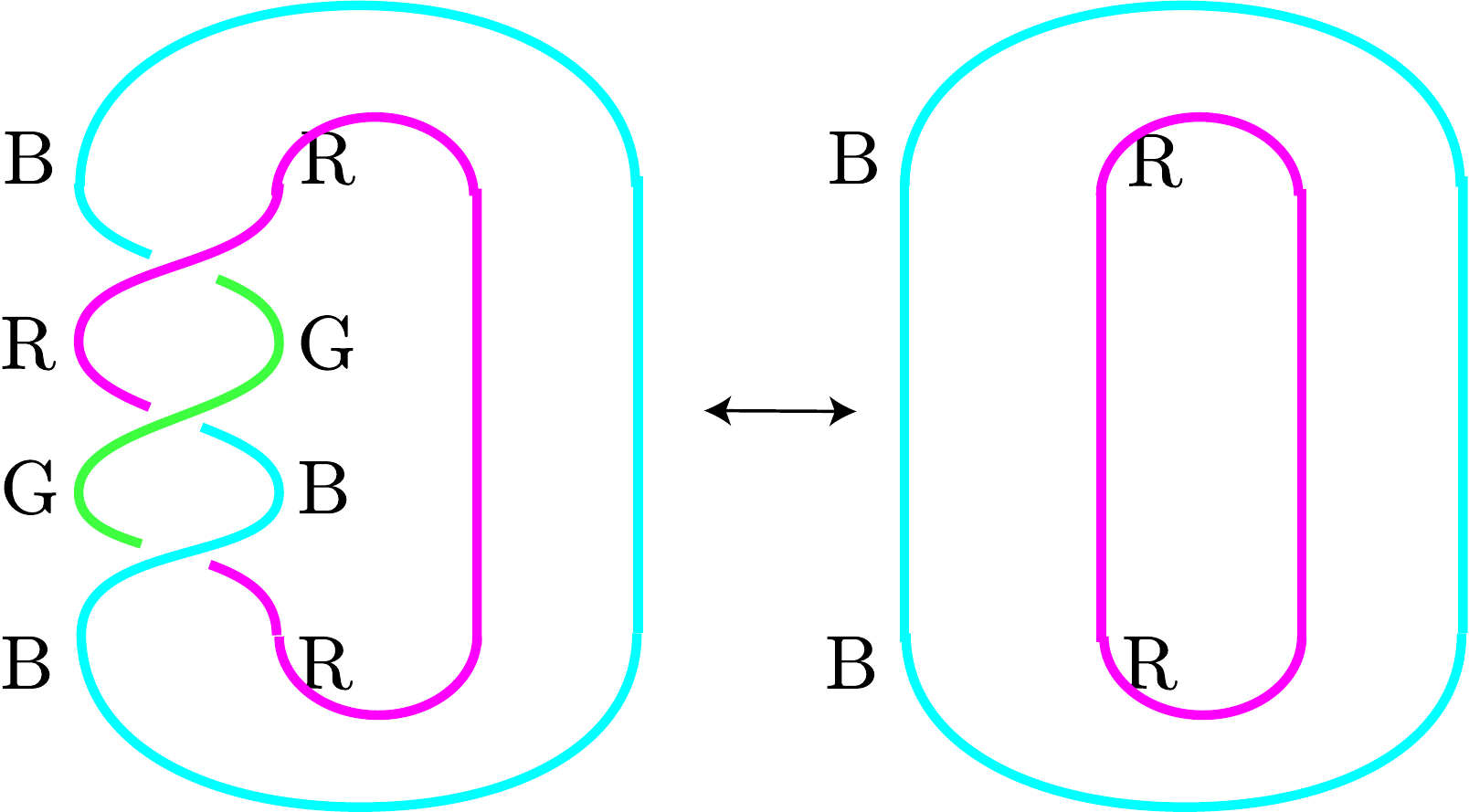}
\end{center}
\caption{Two representations of $S^3$ as a 3-fold simple branched covering}
\label{fig:scovermove_03}
\end{figure}

\begin{definition}\label{def:simplehomoS}{\rm 
A homomorphism $ \rho: \pi_1(S^3 \setminus L, \ast) \to S_d$ sending 
each meridian to a transposition is called a {\it simple} homomorphism.  
}\end{definition}

A link $L$ with a simple homomorphism $ \rho: \pi_1(S^3 \setminus L, \ast) \to S_d$ 
induces a $d$-fold simple branched covering $f : M^3 \to S^3$ branched along $L$. 

Let $D^2 \times S^3$ be a tubular neighborhood of a standardly embedded $S^3$ in $R^5$, and 
let ${\rm pr}: D^2 \times S^3 \to S^3$ be the projection.  

\begin{definition}\label{def:3-dimbraid}{\rm 
A ({\it simple}\/)   ({\it embedded/immersed}\/)  {\it $3$-dimensional braid} is a 
PL map $g : M^3 \to D^2 \times S^3 \subset R^5$ such that 
\begin{itemize}
\item[(1)] 
the composition 
${\rm pr} \circ g: M^3 \to S^3$ is a (simple) branched covering,
\item[(2)] 
$g$ is an embedding/immersion, and 
\item[(3)] 
if $g$ is an immersion, the image of multipoint set under ${\rm pr}$ is a link in $S^3$ avoiding the branch set. 
\end{itemize}
}\end{definition} 

Let $f: M^3 \to S^3$ be a branched covering and  $g : M^3 \to D^2 \times S^3 \subset R^5$ 
an embedded/immersed $3$-dimensional braid.  If  ${\rm pr} \circ g = f$, then we call 
$g$ an {\it embedded/immersed}\/ {\it lift} of $g$.  

\begin{theorem} 
For any $2$-fold simple  branched covering $f: M^3 \to S^3$,  
there exists an embedded lift $g: M^3 \to D^2 \times S^3 \subset R^5$.    
\end{theorem} 

\begin{theorem} 
For any $d$-fold simple  branched covering $f: M^3 \to S^3$,  
there exists an immersed lift $g: M^3 \to D^2 \times S^3 \subset R^5$.    
\end{theorem} 

\begin{problem}{\rm 
When does a simple branched covering $f: M^3 \to S^3$ have an embedded lift? 
}\end{problem}

\underline{In terms of groups} 
\vspace{0.2cm}

Let $L$ be a link in $S^3$. 
Recall Definition~\ref{def:simplehomoS} that a homomorphism $f: \pi_1(S^3 \setminus L) \to S_d$ is {\it simple} if each meridian is mapped to a transposition.  

\begin{definition}\label{def:simplehomoB}{\rm 
A homomorphism $g: \pi_1(S^3 \setminus L) \to B_d$ is {\it simple} if each meridian is mapped to a conjugate of $\sigma_i$ or $\sigma_i^{-1}$.  
}\end{definition}

Let ${\rm pr}: B_d \to S_d$ be the natural projection. 

Let $f: \pi_1(S^3 \setminus L) \to S_d$ and $g : \pi_1(S^3 \setminus L) \to B_d$ be simple homomorphisms. 
If ${\rm pr} \circ g = f$, we say that $g$ is a {\it simple lift} of $f$.  

\begin{problem}{\rm 
Characterize a simple homomorphism $f: \pi_1(S^3 \setminus L) \to S_d$ that has a simple lift.  
}\end{problem} 

\underline{In terms of quandles} 
\vspace{0.2cm}

For an oriented link $L$ in $S^3$, let $Q(S^3, L)$ denote the fundamental quandle of $L$ (\cite{FR1992, Jo1982, Ma1982}).  

Let $T_d$ be the set of transpositions in $S_d$.  
Let $A_d$ be the set of conjugates of standard generators of $B_d$ and their inverses. 
The sets $A_d$ and $T_d$ are regarded as quandles by conjugation.  The natural projection 
${\rm pr}: B_d \to S_d$ induces the projection ${\rm pr}: A_d \to T_d$ which is a surjective quandle homomorphism.  

\begin{problem}{\rm 
Characterize a quandle homomorphism $f: Q(S^3, L) \to T_d$ that has a lift $\tilde f: Q(S^3, L) \to A_d$, 
i.e., ${\rm pr} \circ \tilde f = f$.  
}\end{problem} 

In general we are interested in the following problem.  

\begin{problem}{\rm 
Let $p : \widetilde{Q} \to Q$ be a surjective quandle homomorphism.
Characterize a quandle homomorphism $f : P \to Q$ that
has a lift $\widetilde{f} : P \to \widetilde{Q}$ 
with respect to $p$, {\it i.e.}, $f = p \circ \widetilde{f}$.
}\end{problem}

%%%%%%%%%%%%%%%%%%%%%
\section{$2$-dimensional charts $(m=3)$}
%%%%%%%%%%%%%%%%%%%%%

Permutation charts and braid charts are graphs in $S^2$ describing simple branched coverings of $S^2$ and 
simple 2-dimensional braids.  These notions are generalized into higher dimensions.  The authors are studying 
2-dimensional permutation charts and 2-dimensional braid charts.  They are used to describe 
simple branched coverings of $S^3$ and simple 3-dimensional braids, respectively.

\begin{itemize}
\item A simple embedded branched covering  of $S^3$ $\Leftarrow$ a 2-dimensional  permutation chart.
\item 
A simple  embedded 3-dimensional braid \\ 
$\Leftarrow$ a 2-dimensional  braid chart, or a {\it curtain}.  
\item 
A simple  immersed 3-dimensional braid \\ 
$\Leftarrow$ a 2-dimensional braid chart (or a curtain) with/without {\it nodal curves}.   
\end{itemize}

A $2$-dimensional (permutation or braid) chart  is a 2-dimensional subcomplex of $S^3$ whose faces are (unoriented or oriented), and 
labeled by integers in $\{1, \dots, d-1\}$ such that certain conditions around edges are assumed.  We show some examples of 2-dimensional charts.  

\begin{example}{\rm 
In Figure~\ref{fig:strefoilseif} a trefoil $L$ with a Seifert surface $F$ is depicted.  
When we forget the orientation of $F$, the surface $F$ 
is regarded as a $2$-dimensional permutation chart of degree $2$, or a 
2-dimensional $S_2$-chart.  (We assume that the sheet has label $1$.) 
It induces a monodromy 
$\pi_1(S^3 \setminus L, \ast) \to S_2$ using intersection words.  
It describes a simple embedded $2$-fold branched covering $f_F: M^3 \to S^3$ with branch set $L$.    

When we use the orientation of $F$, the surface $F$ 
is regarded as a $2$-dimensional braid chart of degree $2$, or a 
2-dimensional $B_2$-chart.  (We assume that the sheet has label $1$.) 
It induces a monodromy 
$\pi_1(S^3 \setminus L, \ast) \to B_2$ using intersection words.  
It describes a simple embedded $3$-dimensional braid $g_F:  M^3 \to D^2 \times S^3 \subset R^5$. 
}\end{example}

\begin{figure}[htb]
\begin{center}
\includegraphics[width=2.0in]{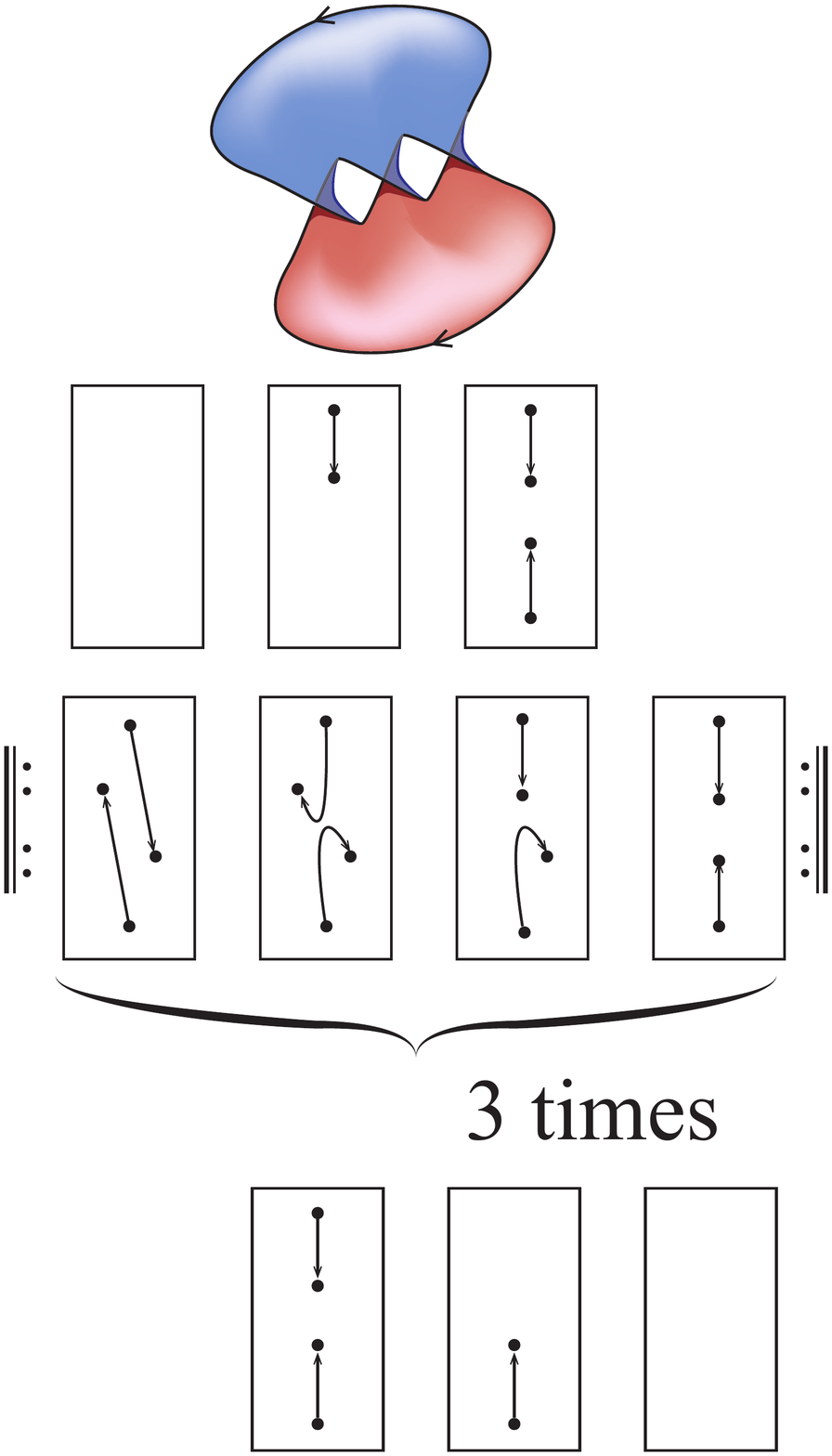}
\end{center}
\caption{A trefoil with a Seifert surface}
\label{fig:strefoilseif}
\end{figure}

\begin{example}{\rm 
In Figure~\ref{fig:sfive2a} a knot $5_2$, denoted by $L$ here, with a Seifert surface, denoted by $F$, is depicted.  
Figure~\ref{fig:sfive2b} shows a motion picture of $L$ and $F$.  

When we forget the orientation of $F$, the surface $F$ 
is regarded as a $2$-dimensional permutation chart of degree $2$, or a 
2-dimensional $S_2$-chart.  (We assume that the sheet has label $1$.) 
It induces a monodromy 
$\pi_1(S^3 \setminus L, \ast) \to S_2$ using intersection words.  
It describes a simple embedded $2$-fold branched covering $f_F: M^3 \to S^3$ with branch set $L$.    

When we use the orientation of $F$, the surface $F$ 
is regarded as a $2$-dimensional braid chart of degree $2$, or a 
2-dimensional $B_2$-chart.  (We assume that the sheet has label $1$.) 
It induces a monodromy 
$\pi_1(S^3 \setminus L, \ast) \to B_2$ using intersection words.  
It describes a simple embedded $3$-dimensional braid $g_F:  M^3 \to D^2 \times S^3 \subset R^5$. 
}\end{example}

\begin{figure}[htb]
\begin{center}
\includegraphics[width=2.2in]{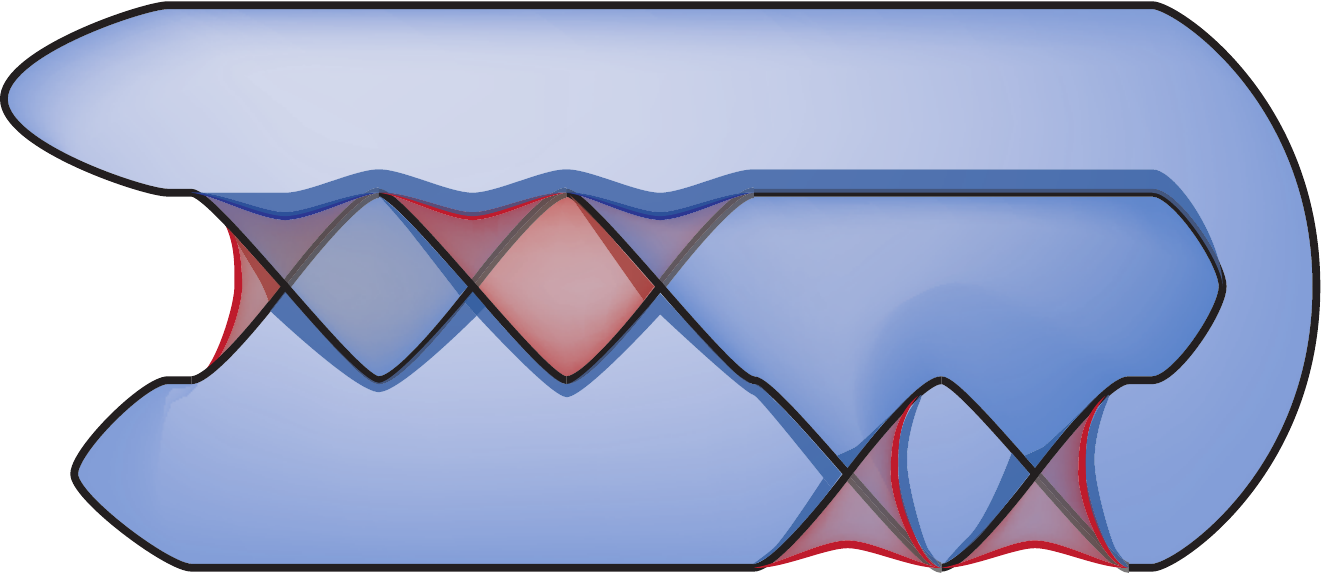}
\end{center}
\caption{A knot $5_2$ with a Seifert surface}
\label{fig:sfive2a}
\end{figure}

\begin{figure}[htb]
\begin{center}
\includegraphics[width=4.5in]{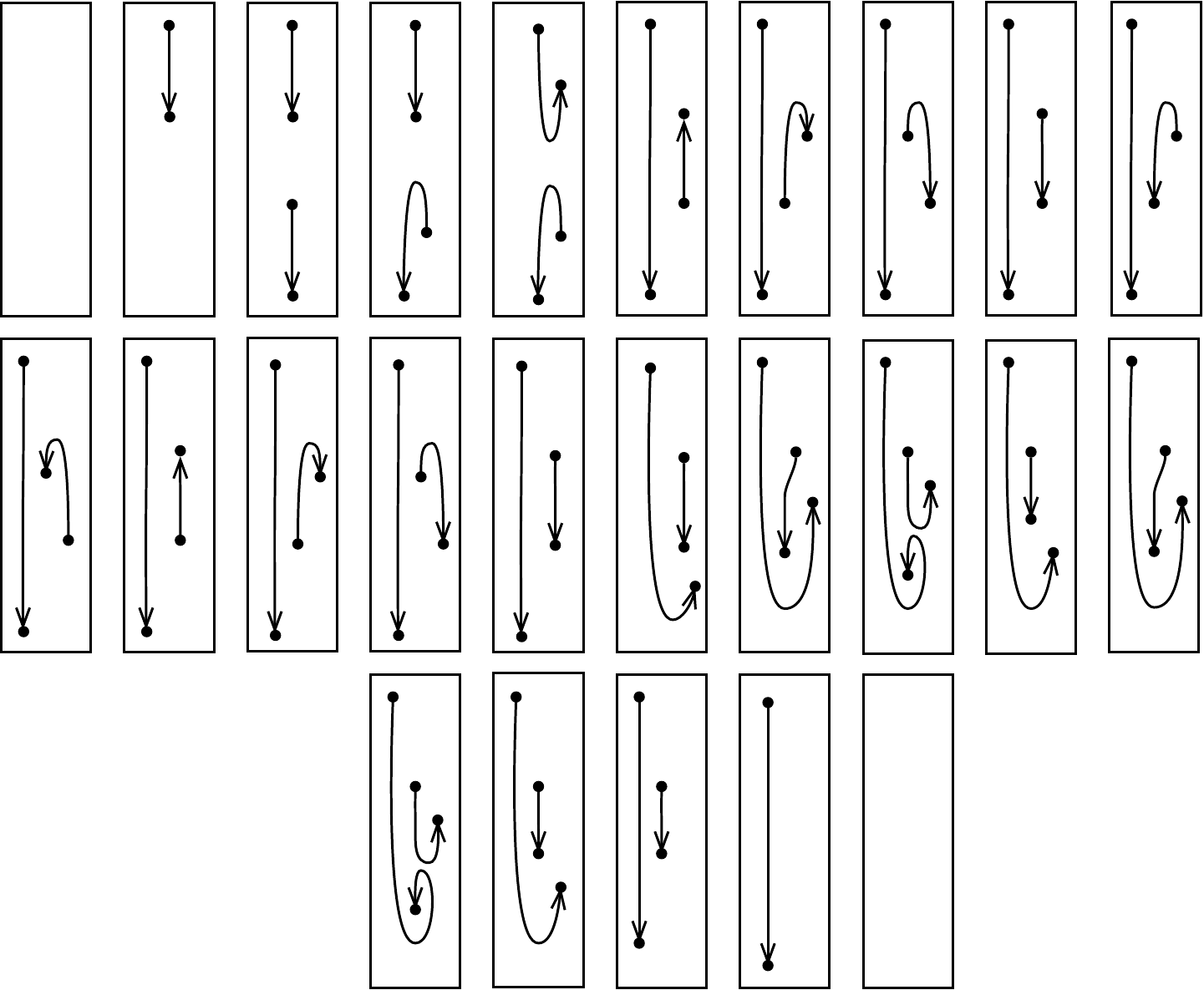}
\end{center}
\caption{A motion picture}
\label{fig:sfive2b}
\end{figure}

\begin{example}{\rm 
Figures~\ref{fig:scurtaintrefoil_01} and \ref{fig:scurtaintrefoil_02} show a $3$-colored trefoil and a $2$-dimensional braid chart. 
Let $L$ be the trefoil knot depicted on the left of Figure~\ref{fig:scurtaintrefoil_01}.  Let $\rho: \pi_1(S^3 \setminus L) \to S_3$ be the monodromy described by the $3$-coloring.   In the right side of Figures~\ref{fig:scurtaintrefoil_01} and \ref{fig:scurtaintrefoil_02}, a motion picture of a $2$-dimensional braid chart $\Gamma$ of degree $3$ is depicted.  The monodromy induced from $\Gamma$ is  $\rho$.  

}\end{example}

\begin{figure}[htb]
\begin{center}
\includegraphics[width=4.in]{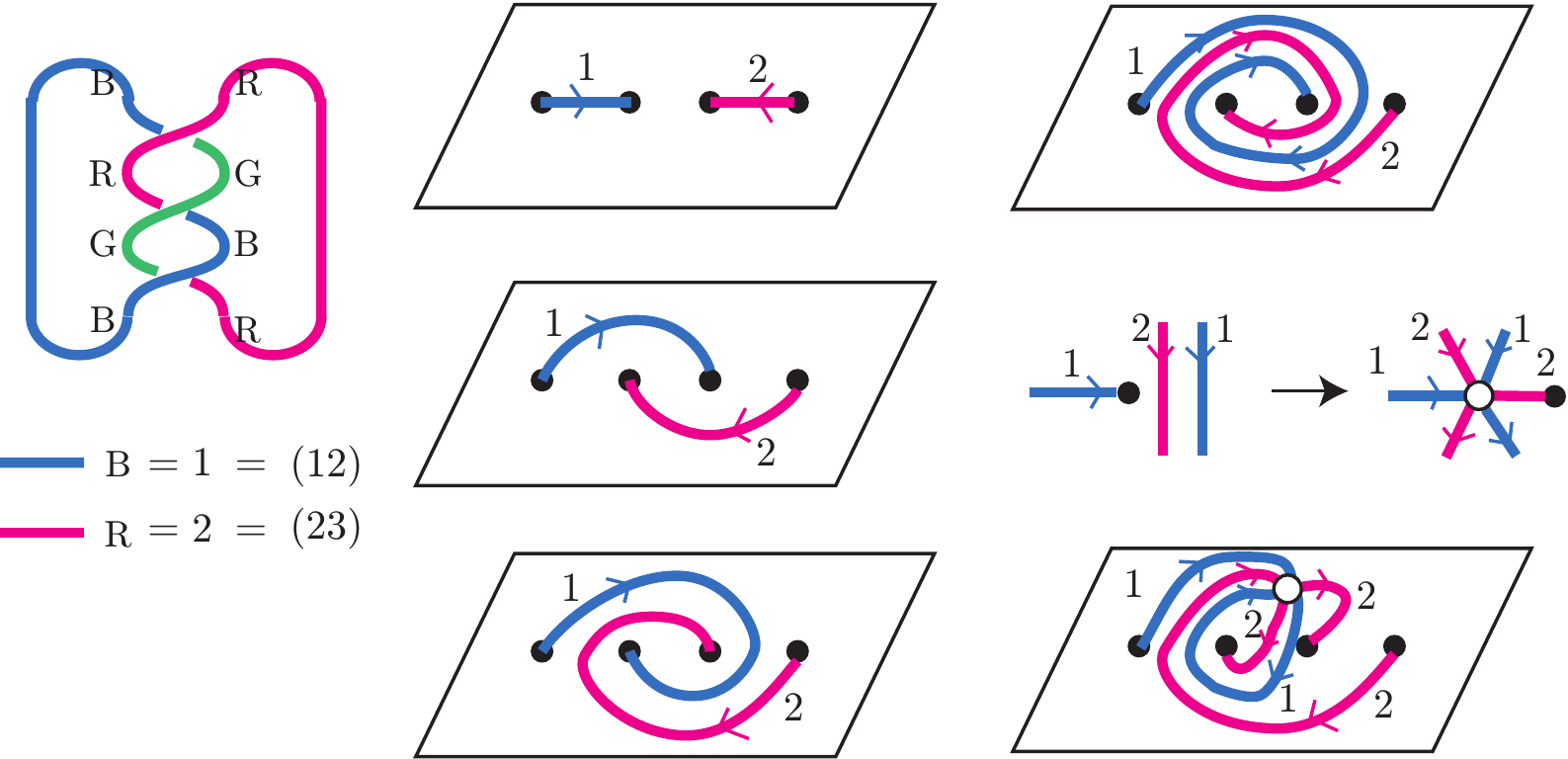}
\end{center}
\caption{A $3$-colored trefoil and a $2$-dimensional braid chart}
\label{fig:scurtaintrefoil_01}
\end{figure}

\begin{figure}[htb]
\begin{center}
\includegraphics[width=4.8in]{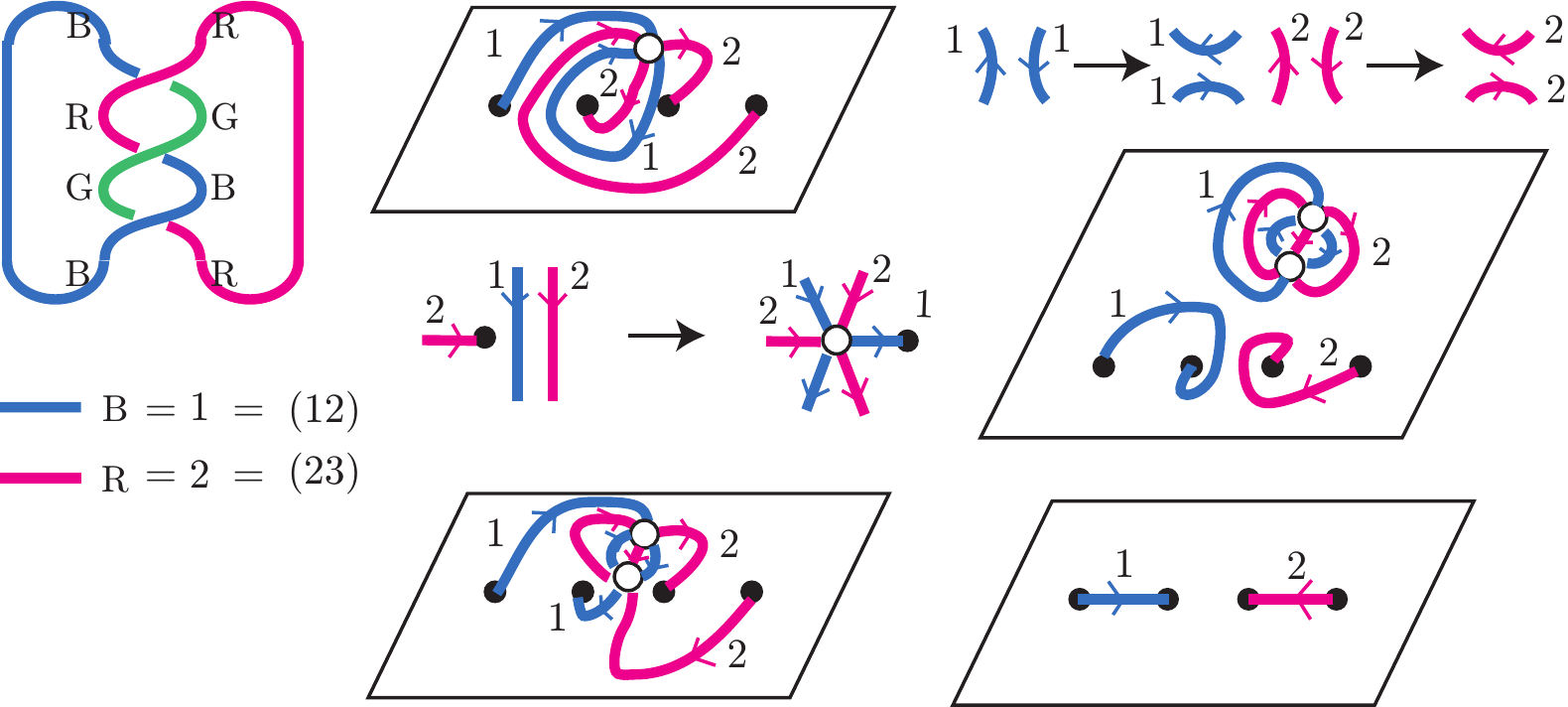}
\end{center}
\caption{A $3$-colored trefoil and a $2$-dimensional braid chart}
\label{fig:scurtaintrefoil_02}
\end{figure}

\subsection*{Acknowledgements.} This  paper was studied with the support of the Ministry of Education Science and Technology (MEST) and the Korean Federation of Science and Technology Societies (KOFST).   SK is being supported by  JSPS grants \#21340015 and \#23654027.

\medskip

\begin{flushleft}
J. Scott Carter \\ 
Department of Mathematics \\ 
University of South Alabama \\ 
Mobile, AL 36688 \\
USA\\ 
E-mail address: {\tt carter@southalabama.edu} 
\end{flushleft}

\begin{flushleft}
Seiichi Kamada \\ 
Department of Mathematics \\ 
Hiroshima University  \\
Hiroshima 739-8526\\ 
JAPAN \\
E-mail address: {\tt  kamada@math.sci.hiroshima-u.ac.jp}
\end{flushleft}

\end{document}